\documentclass[leqno,12pt]{amsart}
\usepackage{amssymb} 
\theoremstyle{plain} 
\newtheorem{theorem}{\indent\sc Theorem}[section] 

\newtheorem{proposition}[theorem]{\indent\sc Proposition}

\theoremstyle{definition} 

\begin{document}

\title[Ordinary differential systems]{Ordinary differential systems in dimension three with affine Weyl group symmetry of types\\
 $D_4^{(1)},B_3^{(1)},G_2^{(1)},D_3^{(2)}$ and $A_2^{(2)}$  \\}
\author{Yusuke Sasano }

\renewcommand{\thefootnote}{\fnsymbol{footnote}}
\footnote[0]{2000\textit{ Mathematics Subjet Classification}.
34M55; 34M45; 58F05; 32S65.}

\maketitle

\begin{abstract}
We present a four-parameter family of ordinary differential systems in dimension three with affine Weyl group symmetry of type $D_4^{(1)}$. By obtaining its first integral, we can reduce this system to the second-order non-linear ordinary differential equations of Painlev\'e type. We also study this system restricted its parameters. Each system can be obtained by connecting some invariant divisors in the system of type $D_4^{(1)}$. Each system admits affine Weyl group symmetry of types $B_3^{(1)},G_2^{(1)},D_3^{(2)}$ and $A_2^{(2)}$, respectively. These symmetries, holomorphy conditions and invariant divisors are new.
\end{abstract}

\section{Introduction}
In this paper, we present a 4-parameter family of ordinary differential systems in dimension three with affine Weyl group symmetry of type $D_4^{(1)}$. By obtaining its first integral, we can reduce this system to the second-order non-linear ordinary differential equations of Painlev\'e type. This reduced system parametrizes the first-order ordinary differential equation:
\begin{equation}
\frac{dX}{dt}=\frac{b(t)}{2\eta}X(X+1)(X+1-\eta)(X-\eta), \quad b(t) \in {\Bbb C}(t), \ \eta \in {\Bbb C}-\{0\}.
\end{equation}

We also study this system restricted its parameters. Each system admits affine Weyl group symmetry of types $B_3^{(1)},G_2^{(1)},D_3^{(2)}$ and $A_2^{(2)}$, respectively.

Each system can be obtained by connecting some invariant divisors in the system of type $D_4^{(1)}$.

The B{\"a}cklund transformations of each system satisfy
\begin{equation}
s_i(g)=g+\frac{\alpha_i}{f_i}\{f_i,g\}+\frac{1}{2!} \left(\frac{\alpha_i}{f_i} \right)^2 \{f_i,\{f_i,g\} \}+\cdots \quad (g \in {\Bbb C}(t)[x,y,z]),
\end{equation}
where poisson bracket $\{,\}$ satisfies the relations:
\begin{equation}
\{z,x\}=\{z,y\}=1, \quad \{x,y\}=0.
\end{equation}
Since these B{\"a}cklund transformations have Lie theoretic origin, similarity reduction of a Drinfeld-Sokolov hierarchy admits such a B{\"a}cklund symmetry.

These symmetries, holomorphy conditions and invariant divisors are new.

\section{The system of type $D_4^{(1)}$}
In this paper, we study the third-order ordinary differential system:
\begin{equation}\label{D4system}
  \left\{
  \begin{aligned}
   \frac{2\eta}{b(t)}\frac{dx}{dt} =&-2 x (y - 1) y z (x - \eta) (2 x - 2 y + 1 - \eta) - 
 2 (2 \alpha_0 + 2 \alpha_1 + 4 \alpha_2 + \alpha_3 + \alpha_4) 
  x^3 y\\
&- 2 (\alpha_0 + \alpha_1) 
  x y^3 + (5 \alpha_0 + 5 \alpha_1 + 
    6 \alpha_2 + \alpha_3 + \alpha_4) 
  x^2 y^2\\
&-\{5 \alpha_0 + 5 \alpha_1 + 6 \alpha_2 + 2 \alpha_3 - 
    3 (2 \alpha_0 + 2 \alpha_1 + 
       4 \alpha_2 + \alpha_3 + \alpha_4) \eta \}x^2 y\\
&+\{3 (\alpha_0 + \alpha_1) - (4 \alpha_0 + 6 \alpha_1 + 
       6 \alpha_2 + \alpha_3 + \alpha_4) \eta \}x y^2 + x^4\\
&+2 (\alpha_0 + \alpha_1 + 2 \alpha_2 + \alpha_3 - \eta) x^3 + 
 2 \alpha_1 \eta y^3\\
&+\{(\alpha_0 + \alpha_1 + 
       2 \alpha_2 + \alpha_3) (\eta^2 - 3 \eta + 
       1) + \alpha_4 \eta^2 \}x^2 + \alpha_1 (\eta - 
    3) \eta y^2\\
&+\{-(\alpha_0 + \alpha_1) + 
    2 (2 \alpha_0 + 3 \alpha_1 + 
       3 \alpha_2 + \alpha_3) \eta - (2 \alpha_0 + 2 \alpha_1 + 
       4 \alpha_2 + \alpha_3 + \alpha_4) \eta^2 \}xy\\
&+(\alpha_0 + \alpha_1 + 2 \alpha_2 + \alpha_3) (\eta - 
    1) \eta x - \alpha_1 (\eta - 1) \eta y,\\
   -\frac{2\eta}{b(t)}\frac{dy}{dt} =&2 x (y - 1) y z (x - \eta) (2 x - 2 y + 1 - \eta) - 
 2 (\alpha_3 + \alpha_4) x^3 y\\
&-2 x y^3 (\alpha_0 + \alpha_1 + 4 \alpha_2 + 2 \alpha_3 + 
    2 \alpha_4) + 
 x^2 y^2 (\alpha_0 + \alpha_1 + 6 \alpha_2 + 5 \alpha_3 + 
    5 \alpha_4)\\
&+x^2 y \{-(\alpha_0 + \alpha_1 + 6 \alpha_2 + 4 \alpha_3 + 
       6 \alpha_4) + 3 (\alpha_3 + \alpha_4) \eta \}\\
&+x y^2\{3 (\alpha_0 + \alpha_1 + 4 \alpha_2 + 2 \alpha_3 + 
       2 \alpha_4) - (2 \alpha_0 + 6 \alpha_2 + 5 \alpha_3 + 
       5 \alpha_4) \eta \}+y^4 + 2 \alpha_4 x^3\\
&+2 y^3 ((\alpha_0 + 2 \alpha_2 + \alpha_3 + \alpha_4) \eta - 
    1) - \alpha_4 (3 \eta - 1) x^2\\
&+y^2\{1 - 3 (\alpha_0 + 
       2 \alpha_2 + \alpha_3 + \alpha_4) \eta + (\alpha_0 + 
       2 \alpha_2 + \alpha_3 + \alpha_4) \eta^2 \}\\
&+x y\{-(\alpha_0 + \alpha_1 + 4 \alpha_2 + 2 \alpha_3 + 
       2 \alpha_4) + 
    2 (\alpha_0 + 3 \alpha_2 + 2 \alpha_3 + 
       3 \alpha_4) \eta - (\alpha_3 + \alpha_4) \eta^2\}\\
&+\alpha_4 (\eta - 1) \eta x - (\alpha_0 + 
    2 \alpha_2 + \alpha_3 + \alpha_4) (\eta - 1) \eta y,\\
   \frac{2\eta}{b(t)}\frac{dz}{dt} =&(2 x - 2 y + 1 - \eta)\{2 x^2 y + 2 x y^2 - x^2 - y^2 \eta - 
    2 (\eta + 1) x y + \eta (x + y) \}z^2\\
&+z[-2 x^3 (\alpha_3 + \alpha_4) + 2 y^3 (\alpha_0 + \alpha_1) + 
    2 x^2 y (\alpha_0 + \alpha_1 + 6 \alpha_2 + 2 \alpha_3 + 
       2 \alpha_4)\\
&-2 x y^2 (2 \alpha_0 + 2 \alpha_1 + 
       6 \alpha_2 + \alpha_3 + \alpha_4) + 
    x^2 \{-(\alpha_0 + \alpha_1 + 6 \alpha_2 + 4 \alpha_3) + 
       3 (\alpha_3 + \alpha_4 ) \eta \}\\
&+y^2 (-3 (\alpha_0 + \alpha_1) + (4 \alpha_0 + 
          6 \alpha_2 + \alpha_3 + \alpha_4) \eta)\\
&-4 x y \{-(\alpha_0 + \alpha_1 + 
          3 \alpha_2 + \alpha_3) + (\alpha_0 + 
          3 \alpha_2 + \alpha_3 + \alpha_4) \eta \}\\
&+x\{-(\alpha_0 + \alpha_1 + 4 \alpha_2 + 2 \alpha_3) + 
       2 (\alpha_0 + 3 \alpha_2 + 
          2 \alpha_3) \eta - (\alpha_3 + \alpha_4) \eta^2 \}\\
&+y\{\alpha_0 + \alpha_1 - 
       2 (2 \alpha_0 + 
          3 \alpha_2 + \alpha_3) \eta + (2 \alpha_0 + 
          4 \alpha_2 + \alpha_3 + \alpha_4) \eta^2 \}\\
&- (\alpha_0+ 2 \alpha_2 + \alpha_3) (\eta - 
       1) \eta] + \alpha_2[(\alpha_0 + \alpha_1 + 
       2 \alpha_2 - \alpha_3 - \alpha_4)x^2\\
&+(\alpha_0 + \alpha_1 - 
       2 \alpha_2 - \alpha_3 - \alpha_4) y^2 - 
    2 (\alpha_0 + \alpha_1 - \alpha_3 - \alpha_4) x y\\
&+x\{\alpha_0 + \alpha_1 - 
       2 \alpha_3 - (2 \alpha_0 + 
          2 \alpha_2 - \alpha_3 - \alpha_4) \eta \}\\
&+y\{-\alpha_0 - \alpha_1 + 2 \alpha_2 + 
       2 \alpha_3 + (2 \alpha_0 - \alpha_3 - \alpha_4) \eta\} +(\eta - 1) (\alpha_2 + \alpha_3 + (\alpha_0 + \alpha_2) \eta)].
   \end{aligned}
  \right. 
\end{equation}
Here $x,y$ and $z$ denote unknown complex variables, $b(t) \in {\Bbb C}(t)$, $\eta \in {\Bbb C}-\{0\}$, and $\alpha_0,\alpha_1, \dots ,\alpha_4$ are complex parameters satisfying the relation:
\begin{equation}
\alpha_0+\alpha_1+2\alpha_2+\alpha_3+\alpha_4=1.
\end{equation}

\begin{theorem}
Let us consider the following ordinary differential system in the polynomial class\rm{:\rm}
\begin{equation*}
  \left\{
  \begin{aligned}
   \frac{dx}{dt} &=f_1(x,y,z),\\
   \frac{dy}{dt} &=f_2(x,y,z),\\
   \frac{dz}{dt} &=f_3(x,y,z).
   \end{aligned}
  \right. 
\end{equation*}
We assume that

$(A1)$ $deg(f_i)=6$ with respect to $x,y,z$.

$(A2)$ The right-hand side of this system becomes again a polynomial in each coordinate system $(x_i,y_i,z_i) \ (i=0,1,\ldots,4)$.
\begin{align}
\begin{split}
0) \ &x_0=x-\eta, \quad y_0=y, \quad z_0=z-\frac{\alpha_0}{x-\eta},\\
1) \ &x_1=x, \quad y_1=y, \quad z_1=z-\frac{\alpha_1}{x},\\
2) \ &x_2=x+\frac{\alpha_2}{z}, \quad y_2=y+\frac{\alpha_2}{z}, \quad z_2=z,\\
3) \ &x_3=x, \quad y_3=y-1, \quad z_3=z-\frac{\alpha_3}{y-1},\\
4) \ &x_4=x, \quad y_4=y, \quad z_4=z-\frac{\alpha_4}{y}.
\end{split}
\end{align}
Then such a system coincides with the system \eqref{D4system}.
\end{theorem}
These transition functions satisfy the condition{\rm:\rm}
\begin{equation*}
dx_i \wedge dy_i \wedge dz_i=dx \wedge dy \wedge dz \quad (i=0,1,\ldots,4).
\end{equation*}

\begin{theorem}
The system \eqref{D4system} admits the affine Weyl group symmetry of type $D_4^{(1)}$ as the group of its B{\"a}cklund transformations, whose generators are explicitly given as follows{\rm : \rm}with the notation $(*):=(x,y,z;\alpha_0,\alpha_1,\ldots,\alpha_4)$,
\begin{align*}
        s_{0}: (*) &\rightarrow \left(x,y,z-\frac{\alpha_0}{x-\eta};-\alpha_0,\alpha_1,\alpha_2+\alpha_0,\alpha_3,\alpha_4 \right),\\
        s_{1}: (*) &\rightarrow \left(x,y,z-\frac{\alpha_1}{x};\alpha_0,-\alpha_1,\alpha_2+\alpha_1,\alpha_3,\alpha_4 \right), \\
        s_{2}: (*) &\rightarrow  \left(x+\frac{\alpha_2}{z},y+\frac{\alpha_2}{z},z;\alpha_0+\alpha_2,\alpha_1+\alpha_2,-\alpha_2,\alpha_3+\alpha_2,\alpha_4+\alpha_2 \right), \\
        s_{3}: (*) &\rightarrow \left(x,y,z-\frac{\alpha_3}{y-1};\alpha_0,\alpha_1,\alpha_2+\alpha_3,-\alpha_3,\alpha_4 \right), \\
        s_{4}: (*) &\rightarrow \left(x,y,z-\frac{\alpha_4}{y};\alpha_0,\alpha_1,\alpha_2+\alpha_4,\alpha_3,-\alpha_4 \right).
\end{align*}
\end{theorem}
We note that the generators $s_0,s_1,\ldots,s_4$ are determined by the invariant divisors \eqref{D4inv} (see next proposition).

\begin{proposition}
The system \eqref{D4system} has the following invariant divisors\rm{:\rm}
\begin{center}\label{D4inv}
\begin{tabular}{|c|c|c|} \hline
parameter's relation & $f_i$ \\ \hline
$\alpha_0=0$ & $f_0:=x-\eta$  \\ \hline
$\alpha_1=0$ & $f_1:=x$  \\ \hline
$\alpha_2=0$ & $f_2:=z$  \\ \hline
$\alpha_3=0$ & $f_3:=y-1$  \\ \hline
$\alpha_4=0$ & $f_4:=y$  \\ \hline
\end{tabular}
\end{center}
\end{proposition}
We note that when $\alpha_1=0$, we see that the system \eqref{D4system} admits a particular solution $x=0$.

\section{Reduction of the $D_4^{(1)}$ system}
In this section, we show that the system \eqref{D4system} has its first integral. Thanks to this first integral, we can reduce the system \eqref{D4system} to the second-order non-linear ordinary differential equations.

\begin{proposition}
The system \eqref{D4system} has its first integral$:$
\begin{equation}\label{eq1}
\frac{d(x-y)}{dt}=\frac{b(t)}{2\eta}(x-y)(x-y+1)(x-y+1-\eta)(x-y-\eta).
\end{equation}
\end{proposition}
Under the condition
\begin{equation}
b(t)=\frac{2\eta}{t(t-1)(t+\eta)(t+\eta-1)},
\end{equation}
the equation \eqref{eq1} admits a particular solution:
\begin{equation}
x=y-t.
\end{equation}

\begin{theorem}
Under the conditions
\begin{align}
&b(t)=\frac{2\eta}{t(t-1)\{t^2+(2\eta-1)t+\eta(\eta-1) \}},\\
&x=y-t,
\end{align}
for the system \eqref{D4system} we make the change of parameters and variables
\begin{equation}\label{65}
\alpha_0=A_1, \quad \alpha_1=A_0, \quad \alpha_2=A_2, \quad \alpha_3=A_3, \quad \alpha_4=A_4, \quad X:=y, \quad Y:=z\\
\end{equation}
from $\alpha_0,\alpha_1,\ldots,\alpha_4,x,y,z$ to $A_i,X,Y$. Then the system \eqref{D4system} can also be written in the new variables $X,Y$ and parameters $A_i$ as a Hamiltonian system. This new system tends to
\begin{align}\label{61}
\begin{split}
\frac{dX}{dt}&=\frac{\partial H_{VI}}{\partial Y}, \quad \frac{dY}{dt}=-\frac{\partial H_{VI}}{\partial X}
\end{split}
\end{align}
with the polynomial Hamiltonian
\begin{align}\label{HVI}
\begin{split}
&H_{VI}(X,Y,t;A_1,A_0,A_2,A_3,A_4)\\
&=\frac{1}{t(t-1)}[Y^2(X-t)(X-1)X-\{(A_1-1)(X-1)X+A_3(X-t)X\\
&+A_4(X-t)(X-1)\}Y+A_2(A_0+A_2)X]
\end{split}
\end{align}
as $\eta \rightarrow \infty$.
\end{theorem}
This system is the Painlev\'e VI system.

We also see that the equation \eqref{eq1} admits a particular solution:
\begin{equation}
x=y.
\end{equation}

\begin{theorem}
Under the condition
\begin{align}
&x=y,
\end{align}
for the system \eqref{D4system} we make the change of variables
\begin{equation}\label{6555}
X:=y, \quad Y:=z\\
\end{equation}
from $x,y,z$ to $X,Y$. Then the system \eqref{D4system} can also be written in the new variables $X,Y$ as a Hamiltonian system. This new system tends to
\begin{align}\label{61555}
\begin{split}
\frac{dX}{dt}&=\frac{\partial H}{\partial Y}, \quad \frac{dY}{dt}=-\frac{\partial H}{\partial X}
\end{split}
\end{align}
with the polynomial Hamiltonian
\begin{align}
\begin{split}
&H(X,Y,t;\alpha_0,\alpha_1,\alpha_2,\alpha_3,\alpha_4)\\
=&-\frac{\eta-1}{2\eta}b(t)[-(X-1)(X-\eta)X^2Y^2-\{2\alpha_2(X-1)(X-\eta)-\alpha_0 \eta(X-1)-\alpha_3(X-\eta)\}XY\\
&-{\alpha_2}^2 X^2+\alpha_2\{1-(\alpha_0+\alpha_1+\alpha_2+\alpha_4)+(1-(\alpha_1+\alpha_2+\alpha_3+\alpha_4))\eta\}X].
\end{split}
\end{align}
\end{theorem}
Elimination of $Y$ from the system \eqref{61555} gives the second-order non-linear ordinary differential equation for the variable $X$:
\begin{align}
\begin{split}
\frac{d^2X}{dt^2}=&\left\{\frac{1}{2(X-1)}+\frac{1}{X}+\frac{1}{2(X-\eta)} \right\} \left(\frac{dX}{dt} \right)^2+\frac{\frac{db(t)}{dt}}{b(t)}\frac{dX}{dt}\\
&-\frac{{b(t)}^2}{8\eta^2(X-1)(X-\eta)}[(\eta-1)^3 X^2\{(\alpha_0 \eta+\alpha_3)X-(1-(\alpha_1+2\alpha_2+\alpha_4))\eta\}\\
& \times \{(\alpha_0 \eta-\alpha_3)X+(1-(2\alpha_0+\alpha_1+2\alpha_2+\alpha_4))\eta \}].
\end{split}
\end{align}

\begin{proposition}
The system \eqref{61555} admits the affine Weyl group symmetry of type $D_4^{(1)}$ as the group of its B{\"a}cklund transformations, whose generators are explicitly given as follows{\rm : \rm}with the notation $(*):=(X,Y;\alpha_0,\alpha_1,\ldots,\alpha_4)$,
\begin{align*}
        s_{0}: (*) &\rightarrow \left(X,Y-\frac{\alpha_0}{X-\eta};-\alpha_0,\alpha_1,\alpha_2+\alpha_0,\alpha_3,\alpha_4 \right),\\
        s_{1}: (*) &\rightarrow \left(X,Y-\frac{\alpha_1}{X};\alpha_0,-\alpha_1,\alpha_2+\alpha_1,\alpha_3,\alpha_4 \right), \\
        s_{2}: (*) &\rightarrow  \left(X+\frac{\alpha_2}{Y},Y;\alpha_0+\alpha_2,\alpha_1+\alpha_2,-\alpha_2,\alpha_3+\alpha_2,\alpha_4+\alpha_2 \right), \\
        s_{3}: (*) &\rightarrow \left(X,Y-\frac{\alpha_3}{X-1};\alpha_0,\alpha_1,\alpha_2+\alpha_3,-\alpha_3,\alpha_4 \right), \\
        s_{4}: (*) &\rightarrow \left(X,Y-\frac{\alpha_4}{X};\alpha_0,\alpha_1,\alpha_2+\alpha_4,\alpha_3,-\alpha_4 \right).
\end{align*}
\end{proposition}

\section{The system of type $B_3^{(1)}$}
In this section, we present a 3-parameter family of ordinary differential systems in dimension three with affine Weyl group symmetry of type $B_3^{(1)}$. This system is equivalent to the system \eqref{D4system} restricted its parameters.

This system can be obtained by connecting the invariant divisors $x$ and $y$ in the system \eqref{D4system}.

For this system, we can discuss some reductions to the second-order ordinary differential equations in the same way in previous section because this system is equivalent to the system \eqref{D4system} only restricted its parameters.

\begin{theorem}
Let us consider the following ordinary differential system in the polynomial class\rm{:\rm}
\begin{equation*}
  \left\{
  \begin{aligned}
   \frac{dx}{dt} &=f_1(x,y,z),\\
   \frac{dy}{dt} &=f_2(x,y,z),\\
   \frac{dz}{dt} &=f_3(x,y,z).
   \end{aligned}
  \right. 
\end{equation*}
We assume that

$(A1)$ $deg(f_i)=6$ with respect to $x,y,z$.

$(A2)$ The right-hand side of this system becomes again a polynomial in each coordinate system $(x_i,y_i,z_i) \ (i=0,1,2,3)$.
\begin{align}
\begin{split}
0) \ &x_0=x-\eta, \quad y_0=y, \quad z_0=z-\frac{\beta_0}{x-\eta},\\
1) \ &x_1=x, \quad y_1=y-1, \quad z_1=z-\frac{\beta_1}{y-1},\\
2) \ &x_2=x+\frac{\beta_2}{z}, \quad y_2=y+\frac{\beta_2}{z}, \quad z_2=z,\\
3) \ &x_3=x, \quad y_3=y, \quad z_3=z-\frac{\beta_3(x+y)}{xy}.
\end{split}
\end{align}
Then such a system coincides with the system \eqref{D4system} with the parameter's relations:
\begin{equation}
\alpha_4=\alpha_1, \quad \beta_0:=\alpha_0, \quad \beta_1:=\alpha_3, \quad \beta_2:=\alpha_2, \quad \beta_3:=\alpha_1.
\end{equation}
\end{theorem}
Here, the complex parameters $\beta_0,\beta_1,\beta_2,\beta_3$ satisfy the relation:
\begin{equation}
\beta_0+\beta_1+2\beta_2+2\beta_3=1.
\end{equation}

These transition functions satisfy the condition{\rm:\rm}
\begin{equation*}
dx_i \wedge dy_i \wedge dz_i=dx \wedge dy \wedge dz \quad (i=0,1,2,3).
\end{equation*}

\begin{theorem}
This system admits the affine Weyl group symmetry of type $B_3^{(1)}$ as the group of its B{\"a}cklund transformations, whose generators are explicitly given as follows{\rm : \rm}with the notation $(*):=(x,y,z;\beta_0,\beta_1,\beta_2,\beta_3)$,
\begin{align*}
        s_{0}: (*) &\rightarrow \left(x,y,z-\frac{\beta_0}{x-\eta};-\beta_0,\beta_1,\beta_2+\beta_0,\beta_3 \right),\\
        s_{1}: (*) &\rightarrow \left(x,y,z-\frac{\beta_1}{y-1};\beta_0,-\beta_1,\beta_2+\beta_1,\beta_3 \right), \\
        s_{2}: (*) &\rightarrow  \left(x+\frac{\beta_2}{z},y+\frac{\beta_2}{z},z;\beta_0+\beta_2,\beta_1+\beta_2,-\beta_2,\beta_3+\beta_2 \right), \\
        s_{3}: (*) &\rightarrow \left(x,y,z-\frac{\beta_3(x+y)}{xy};\beta_0,\beta_1,\beta_2+2\beta_3,-\beta_3 \right).
\end{align*}
\end{theorem}
The B{\"a}cklund transformations of each system satisfy
\begin{equation}
s_i(g)=g+\frac{\alpha_i}{f_i}\{f_i,g\}+\frac{1}{2!} \left(\frac{\alpha_i}{f_i} \right)^2 \{f_i,\{f_i,g\} \}+\cdots \quad (g \in {\Bbb C}(t)[x,y,z]),
\end{equation}
where poisson bracket $\{,\}$ satisfies the relations:
\begin{equation}
\{z,x\}=\{z,y\}=1, \quad \{x,y\}=0.
\end{equation}
Since these B{\"a}cklund transformations have Lie theoretic origin, similarity reduction of a Drinfeld-Sokolov hierarchy admits such a B{\"a}cklund symmetry.

We note that all generators $s_0,s_1,s_2,s_3$ are determined by the invariant divisors \eqref{B3inv} (see next proposition).

\begin{proposition}
This system has the following invariant divisors\rm{:\rm}
\begin{center}\label{B3inv}
\begin{tabular}{|c|c|c|} \hline
parameter's relation & $f_i$ \\ \hline
$\beta_0=0$ & $f_0:=x-\eta$  \\ \hline
$\beta_1=0$ & $f_1:=y-1$  \\ \hline
$\beta_2=0$ & $f_2:=z$  \\ \hline
$\beta_3=0$ & $f_3:=xy$  \\ \hline
\end{tabular}
\end{center}
\end{proposition}
We note that when $\beta_3=0$, after we make the birational transformations:
\begin{equation}
x_3=xy, \ y_3=y, \ z_3=z
\end{equation}
we see that in the coordinate system $(x_3,y_3,z_3)$ the system admits a particular solution $x_3=0$.

\section{The system of type $D_3^{(2)}$}
In this section, we present a 2-parameter family of ordinary differential systems in dimension three with affine Weyl group symmetry of type $D_3^{(2)}$. This system is equivalent to the system \eqref{D4system} restricted its parameters.

\begin{theorem}
Let us consider the following ordinary differential system in the polynomial class\rm{:\rm}
\begin{equation*}
  \left\{
  \begin{aligned}
   \frac{dx}{dt} &=f_1(x,y,z),\\
   \frac{dy}{dt} &=f_2(x,y,z),\\
   \frac{dz}{dt} &=f_3(x,y,z).
   \end{aligned}
  \right. 
\end{equation*}
We assume that

$(A1)$ $deg(f_i)=6$ with respect to $x,y,z$.

$(A2)$ The right-hand side of this system becomes again a polynomial in each coordinate system $(x_i,y_i,z_i) \ (i=1,2,3)$.
\begin{align}
\begin{split}
1) \ &x_1=x-\eta, \quad y_1=y-1, \quad z_1=z-\frac{\beta_0(x+y-\eta-1)}{(x-\eta)(y-1)},\\
2) \ &x_2=x+\frac{\beta_1}{z}, \quad y_2=y+\frac{\beta_1}{z}, \quad z_2=z,\\
3) \ &x_3=x, \quad y_3=y, \quad z_3=z-\frac{\beta_2(x+y)}{xy}.
\end{split}
\end{align}
Then such a system coincides with the system \eqref{D4system} with the parameter's relations:
\begin{equation}
\alpha_4=\alpha_1, \quad \alpha_0=\alpha_3, \quad \beta_0:=\alpha_1, \quad \beta_1:=\alpha_2, \quad \beta_2:=\alpha_3.
\end{equation}
\end{theorem}
Here, the complex parameters $\beta_0,\beta_1,\beta_2$ satisfy the relation:
\begin{equation}
\beta_0+\beta_1+\beta_2=\frac{1}{2}.
\end{equation}

These transition functions satisfy the condition{\rm:\rm}
\begin{equation*}
dx_i \wedge dy_i \wedge dz_i=dx \wedge dy \wedge dz \quad (i=1,2,3).
\end{equation*}

\begin{theorem}
This system admits the affine Weyl group symmetry of type $D_3^{(2)}$ as the group of its B{\"a}cklund transformations, whose generators are explicitly given as follows{\rm : \rm}with the notation $(*):=(x,y,z;\beta_0,\beta_1,\beta_2)$,
\begin{align*}
        s_{0}: (*) &\rightarrow \left(x,y,z-\frac{\beta_0(x+y-\eta-1)}{(x-\eta)(y-1)};-\beta_0,\beta_1+2\beta_0,\beta_2 \right),\\
        s_{1}: (*) &\rightarrow \left(x+\frac{\beta_1}{z},y+\frac{\beta_1}{z},z;\beta_0+\beta_1,-\beta_1,\beta_2+\beta_1 \right), \\
        s_{2}: (*) &\rightarrow  \left(x,y,z-\frac{\beta_2(x+y)}{xy};\beta_0,\beta_1+2\beta_2,-\beta_2 \right).
\end{align*}
\end{theorem}
We note that the generators $s_0,s_1,s_2$ are determined by the invariant divisors \eqref{D3inv} (see next proposition).

\begin{proposition}
This system has the following invariant divisors\rm{:\rm}
\begin{center}\label{D3inv}
\begin{tabular}{|c|c|c|} \hline
parameter's relation & $f_i$ \\ \hline
$\beta_0=0$ & $f_0:=(x-\eta)(y-1)$  \\ \hline
$\beta_1=0$ & $f_1:=z$  \\ \hline
$\beta_2=0$ & $f_2:=xy$  \\ \hline
\end{tabular}
\end{center}
\end{proposition}

\section{The system of type $G_2^{(1)}$}
In this section, we present a 2-parameter family of ordinary differential systems in dimension three with affine Weyl group symmetry of type $G_2^{(1)}$. This system is equivalent to the system \eqref{D4system} restricted its parameters
\begin{equation}
\alpha_4=\alpha_3=\alpha_1, \quad \beta_0:=\alpha_0, \quad \beta_1:=\alpha_2, \quad \beta_2:=\alpha_3.
\end{equation}
Here, the complex parameters $\beta_0,\beta_1,\beta_2$ satisfy the relation:
\begin{equation}
\beta_0+2\beta_1+3\beta_2=1.
\end{equation}

\begin{theorem}
This system admits the affine Weyl group symmetry of type $G_2^{(1)}$ as the group of its B{\"a}cklund transformations, whose generators are explicitly given as follows{\rm : \rm}with the notation $(*):=(x,y,z;\beta_0,\beta_1,\beta_2)$,
\begin{align*}
        s_{0}: (*) &\rightarrow \left(x,y,z-\frac{\beta_0}{x-\eta};-\beta_0,\beta_1+\beta_0,\beta_2 \right),\\
        s_{1}: (*) &\rightarrow \left(x+\frac{\beta_1}{z},y+\frac{\beta_1}{z},z;\beta_0+\beta_1,-\beta_1,\beta_2+\beta_1 \right), \\
        s_{2}: (*) &\rightarrow  \left(x,y,z-\frac{\beta_2 \{y(y-1)+x(y-1)+xy \}}{xy(y-1)};\beta_0,\beta_1+3\beta_2,-\beta_2 \right).
\end{align*}
\end{theorem}
We note that the generators $s_0,s_1,s_2$ are determined by the invariant divisors \eqref{G2inv} (see next proposition).

\begin{proposition}
This system has the following invariant divisors\rm{:\rm}
\begin{center}\label{G2inv}
\begin{tabular}{|c|c|c|} \hline
parameter's relation & $f_i$ \\ \hline
$\beta_0=0$ & $f_0:=x-\eta$  \\ \hline
$\beta_1=0$ & $f_1:=z$  \\ \hline
$\beta_2=0$ & $f_2:=xy(y-1)$  \\ \hline
\end{tabular}
\end{center}
\end{proposition}

\section{The system of type $A_2^{(2)}$}
In this section, we present a 1-parameter family of ordinary differential systems in dimension three with affine Weyl group symmetry of type $A_2^{(2)}$. This system is equivalent to the system \eqref{D4system} restricted its parameters
\begin{equation}
\alpha_4=\alpha_3=\alpha_0=\alpha_1, \quad \beta_0:=\alpha_2, \quad \beta_1:=\alpha_1.
\end{equation}
Here, the complex parameters $\beta_0,\beta_1$ satisfy the relation:
\begin{equation}
\beta_0+2\beta_1=\frac{1}{2}.
\end{equation}

\begin{theorem}
This system admits the affine Weyl group symmetry of type $A_2^{(2)}$ as the group of its B{\"a}cklund transformations, whose generators are explicitly given as follows{\rm : \rm}with the notation $(*):=(x,y,z;\beta_0,\beta_1)$,
\begin{align*}
        s_{0}: (*) \rightarrow& \left(x+\frac{\beta_0}{z},y+\frac{\beta_0}{z},z;-\beta_0,\beta_1+\beta_0 \right),\\
        s_{1}: (*) \rightarrow& (x,y,z-\frac{\beta_1\{y(x-\eta)(y-1)+x(x-\eta)(y-1)+xy(y-1)+xy(x-\eta)\}}{xy(x-\eta)(y-1)};\\
&\beta_0+4\beta_1,-\beta_1).
\end{align*}
\end{theorem}
We note that the generators $s_0,s_1$ are determined by the invariant divisors \eqref{A2inv} (see next proposition).

\begin{proposition}
This system has the following invariant divisors\rm{:\rm}
\begin{center}\label{A2inv}
\begin{tabular}{|c|c|c|} \hline
parameter's relation & $f_i$ \\ \hline
$\beta_0=0$ & $f_0:=z$  \\ \hline
$\beta_1=0$ & $f_1:=xy(x-\eta)(y-1)$  \\ \hline
\end{tabular}
\end{center}
\end{proposition}

\section{Appendix A} 
It is well-known that the fifth Painlev\'e equation has symmetries under the affine Weyl group of type $A_3^{(1)}$. In this section, we present a 3-parameter family of the systems of the first-order ordinary differential equations.
\begin{theorem}
The fifth Painlev\'e equation is equivalent to a 3-parameter family of the systems of the first-order ordinary differential equations$:$
\begin{equation}\label{symmetricPV}
  \left\{
  \begin{aligned}
   \frac{df_0}{dt} &=-\frac{\varphi}{2}(-2f_0f_1f_2+af_0f_2+(\alpha_0+\alpha_1+\alpha_3)f_0-\alpha_0f_2),\\
   \frac{df_1}{dt} &=\frac{\varphi}{2}(-f_0f_1^2-f_1^2f_2+af_0f_1+af_1f_2-a\alpha_1+(\alpha_1+\alpha_3)f_1),\\
   \frac{df_2}{dt} &=-\frac{\varphi}{2}(-2f_0f_1f_2+af_0f_2+(\alpha_1+\alpha_2+\alpha_3)f_2-\alpha_2f_0).
   \end{aligned}
  \right. 
\end{equation}
Here, $f_0,f_1$ and $f_2$ denote unknown complex variables and $\alpha_0,\alpha_1,\alpha_2,\alpha_3$ and $a$ are constant complex parameters with $\alpha_0+\alpha_1+\alpha_2+\alpha_3=1$ and $\varphi$ is a nonzero parameter which can be fixed arbitrarily.
\end{theorem}
We see that the system \eqref{symmetricPV} has its first integral:
\begin{equation}
\frac{d(f_2-f_0)}{dt}=f_2-f_0.
\end{equation}
We can solve this equation by
\begin{equation}
f_2-f_0=e^{(t+c)}.
\end{equation}
Here we set
\begin{equation}
t+c=logT, \ x:=f_0, \ y:=f_1,
\end{equation}
then we can obtain the fifth Painlev\'e system:
\begin{equation}
  \left\{
  \begin{aligned}
   \frac{dx}{dT} &=\frac{\partial H_{V}}{\partial y}=-\frac{2x^2y}{T}+\frac{ax^2}{T}-2xy+\left(a+\frac{\alpha_1+\alpha_3}{T} \right)x-\alpha_0,\\
   \frac{dy}{dT} &=-\frac{\partial H_{V}}{\partial x}=\frac{2xy^2}{T}+y^2-\frac{2a xy}{T}-\left(a+\frac{\alpha_1+\alpha_3}{T} \right)y+\frac{a\alpha_1}{T}
   \end{aligned}
  \right. 
\end{equation}
with the polynomial Hamiltonian
\begin{equation}
H_{V}=-\frac{x^2y^2}{T}+\frac{a x^2y}{T}-xy^2+\left(a+\frac{\alpha_1+\alpha_3}{T} \right)xy-\alpha_0 y-\frac{a\alpha_1}{T}x.
\end{equation}

\begin{theorem}
This system admits the affine Weyl group symmetry of type $A_3^{(1)}$ as the group of its B{\"a}cklund transformations, whose generators $s_0,s_1,s_2,s_3,\pi$ are explicitly given as follows{\rm : \rm}
\begin{align}
\begin{split}
s_0: &(f_0,f_1,f_2;\alpha_0,\alpha_1,\alpha_2,\alpha_3) \rightarrow \left(f_0,f_1+\frac{\alpha_0}{f_0},f_2;-\alpha_0,\alpha_1+\alpha_0,\alpha_2,\alpha_3+\alpha_0 \right),\\
s_1: &(f_0,f_1,f_2;\alpha_0,\alpha_1,\alpha_2,\alpha_3) \rightarrow \left(f_0-\frac{\alpha_1}{f_1},f_1,f_2-\frac{\alpha_1}{f_1};\alpha_0+\alpha_1,-\alpha_1,\alpha_2+\alpha_1,\alpha_3 \right),\\
s_2: &(f_0,f_1,f_2;\alpha_0,\alpha_1,\alpha_2,\alpha_3) \rightarrow \left(f_0,f_1+\frac{\alpha_2}{f_2},f_2;\alpha_0,\alpha_1+\alpha_2,-\alpha_2,\alpha_3+\alpha_2 \right),\\
s_3: &(f_0,f_1,f_2;\alpha_0,\alpha_1,\alpha_2,\alpha_3) \rightarrow \left(f_0-\frac{\alpha_3}{f_1-a},f_1,f_2-\frac{\alpha_3}{f_1-a};\alpha_0+\alpha_3,\alpha_1,\alpha_2+\alpha_3,-\alpha_3 \right),\\
\pi: &(f_0,f_1,f_2;\alpha_0,\alpha_1,\alpha_2,\alpha_3) \rightarrow (f_2,f_1,f_0;\eta,\alpha_2,\alpha_1,\alpha_0,\alpha_3).
\end{split}
\end{align}
\end{theorem}
The B{\"a}cklund transformations of each system satisfy
\begin{equation}
s_i(g)=g+\frac{\alpha_i}{f_i}\{f_i,g\}+\frac{1}{2!} \left(\frac{\alpha_i}{f_i} \right)^2 \{f_i,\{f_i,g\} \}+\cdots \quad (g \in {\Bbb C}[f_0,f_1,f_2]),
\end{equation}
where poisson bracket $\{,\}$ satisfies the relations:
\begin{equation}
\{f_0,f_1\}=\{f_2,f_1\}=1, \quad \{f_0,f_2\}=0.
\end{equation}
Since these B{\"a}cklund transformations have Lie theoretic origin, similarity reduction of a Drinfeld-Sokolov hierarchy admits such a B{\"a}cklund symmetry.

\section{Appendix B}
It is well-known that the third Painlev\'e equation has symmetries under the affine Weyl group of type $C_2^{(1)}$. In this section, we present a new representation of the third Painlev\'e equation.
\begin{theorem}
The third Painlev\'e equation can be written in the following symmetric form$:$
\begin{equation}\label{symmetricPIII}
  \left\{
  \begin{aligned}
   \frac{df_0}{dt} &=-2f_0f_1f_2+(\alpha_0+2\alpha_1)f_0-\alpha_0f_2,\\
   \frac{df_1}{dt} &=(f_0+f_2)f_1^2-2\alpha_1f_1+\eta,\\
   \frac{df_2}{dt} &=-2f_0f_1f_2+(2\alpha_1+\alpha_2)f_2-\alpha_2f_0.
   \end{aligned}
  \right. 
\end{equation}
Here, $f_0,f_1$ and $f_2$ denote unknown complex variables and $\alpha_0,\alpha_1$ and $\alpha_2$ are constant parameters with $\alpha_0+2\alpha_1+\alpha_2=1$ and $\eta$ is a nonzero parameter which can be fixed arbitrarily.
\end{theorem}
We see that the system \eqref{symmetricPIII} has its first integral:
\begin{equation}
\frac{d(f_0-f_2)}{dt}=f_0-f_2.
\end{equation}
We can solve this equation by
\begin{equation}
f_0-f_2=e^{(t+c)}.
\end{equation}
Here we set
\begin{equation}
t+c=logT, \ x:=\frac{1}{f_1}, \ y:=-(f_1f_2+\alpha_2)f_1,
\end{equation}
then we can obtain the third Painlev\'e system:
\begin{equation}
  \left\{
  \begin{aligned}
   \frac{dx}{dT} &=\frac{\partial H_{III}}{\partial y}=\frac{2x^2y}{T}-\frac{\eta x^2}{T}+\frac{2(\alpha_1+\alpha_2)x}{T}-1,\\
   \frac{dy}{dT} &=-\frac{\partial H_{III}}{\partial x}=-\frac{2xy^2}{T}+\frac{2\eta xy}{T}-\frac{2(\alpha_1+\alpha_2)y}{T}+\frac{\eta \alpha_2}{T}
   \end{aligned}
  \right. 
\end{equation}
with the polynomial Hamiltonian
\begin{equation}
H_{III}=\frac{x^2y^2-\eta x^2y+2(\alpha_1+\alpha_2)xy-Ty-\eta \alpha_2 x}{T}.
\end{equation}

\begin{theorem}\label{th:1}
This system admits extended affine Weyl group symmetry of type $C_2^{(1)}$ as the group of its B{\"a}cklund transformations, whose generators $s_0,s_1,s_2,\pi$ are explicitly given as follows{\rm : \rm}
\begin{align}
\begin{split}
s_0: (f_0,f_1,f_2;\eta,\alpha_0,\alpha_1,\alpha_2) \rightarrow &\left(f_0,f_1+\frac{\alpha_0}{f_0},f_2;\eta,-\alpha_0,\alpha_1+\alpha_0,\alpha_2 \right),\\
s_1: (f_0,f_1,f_2;\eta,\alpha_0,\alpha_1,\alpha_2) \rightarrow &(f_0-\frac{2\alpha_1}{f_1}+\frac{\eta}{f_1^2},f_1,f_2-\frac{2\alpha_1}{f_1}+\frac{\eta}{f_1^2};\\
&-\eta,\alpha_0+2\alpha_1,-\alpha_1,\alpha_2+2\alpha_1),\\
s_2: (f_0,f_1,f_2;\eta,\alpha_0,\alpha_1,\alpha_2) \rightarrow &\left(f_0,f_1+\frac{\alpha_2}{f_2},f_2;\eta,\alpha_0,\alpha_1+\alpha_2,-\alpha_2 \right),\\
\pi: (f_0,f_1,f_2;\eta,\alpha_0,\alpha_1,\alpha_2) \rightarrow &(f_2,f_1,f_0;\eta,\alpha_2,\alpha_1,\alpha_0).
\end{split}
\end{align}
\end{theorem}

\begin{figure}[h]
\unitlength 0.1in
\begin{picture}(29.33,16.80)(9.50,-17.50)
%
\special{pn 20}%
\special{ar 1370 1060 263 263  0.0000000 6.2831853}%
%
\special{pn 20}%
\special{ar 2490 1070 263 263  0.0000000 6.2831853}%
%
\special{pn 20}%
\special{ar 3620 1060 263 263  0.0000000 6.2831853}%
%
\special{pn 20}%
\special{pa 1620 980}%
\special{pa 2210 980}%
\special{fp}%
\special{sh 1}%
\special{pa 2210 980}%
\special{pa 2143 960}%
\special{pa 2157 980}%
\special{pa 2143 1000}%
\special{pa 2210 980}%
\special{fp}%
%
\special{pn 20}%
\special{pa 1620 1190}%
\special{pa 2230 1190}%
\special{fp}%
\special{sh 1}%
\special{pa 2230 1190}%
\special{pa 2163 1170}%
\special{pa 2177 1190}%
\special{pa 2163 1210}%
\special{pa 2230 1190}%
\special{fp}%
%
\special{pn 20}%
\special{pa 3370 990}%
\special{pa 2770 990}%
\special{fp}%
\special{sh 1}%
\special{pa 2770 990}%
\special{pa 2837 1010}%
\special{pa 2823 990}%
\special{pa 2837 970}%
\special{pa 2770 990}%
\special{fp}%
%
\special{pn 20}%
\special{pa 3400 1210}%
\special{pa 2760 1210}%
\special{fp}%
\special{sh 1}%
\special{pa 2760 1210}%
\special{pa 2827 1230}%
\special{pa 2813 1210}%
\special{pa 2827 1190}%
\special{pa 2760 1210}%
\special{fp}%
\put(12.3000,-11.7000){\makebox(0,0)[lb]{$f_0$}}%
\put(34.9000,-11.6000){\makebox(0,0)[lb]{$f_2$}}%
\put(9.5000,-8.0000){\makebox(0,0)[lb]{$\alpha_0$}}%
\put(20.8000,-8.1000){\makebox(0,0)[lb]{$\alpha_1$}}%
\put(32.2000,-8.2000){\makebox(0,0)[lb]{$\alpha_2$}}%
%
\special{pn 8}%
\special{pa 2460 260}%
\special{pa 2460 1750}%
\special{dt 0.045}%
\special{pa 2460 1750}%
\special{pa 2460 1749}%
\special{dt 0.045}%
%
\special{pn 20}%
\special{pa 2410 510}%
\special{pa 2367 505}%
\special{pa 2326 499}%
\special{pa 2289 492}%
\special{pa 2257 482}%
\special{pa 2232 470}%
\special{pa 2215 453}%
\special{pa 2209 431}%
\special{pa 2214 405}%
\special{pa 2228 377}%
\special{pa 2251 349}%
\special{pa 2281 323}%
\special{pa 2316 301}%
\special{pa 2356 284}%
\special{pa 2398 273}%
\special{pa 2442 266}%
\special{pa 2487 264}%
\special{pa 2529 265}%
\special{pa 2570 270}%
\special{pa 2606 279}%
\special{pa 2636 290}%
\special{pa 2660 305}%
\special{pa 2675 321}%
\special{pa 2681 340}%
\special{pa 2676 360}%
\special{pa 2660 382}%
\special{pa 2637 405}%
\special{pa 2609 428}%
\special{pa 2580 450}%
\special{sp}%
%
\special{pn 20}%
\special{pa 2650 380}%
\special{pa 2520 490}%
\special{fp}%
\special{sh 1}%
\special{pa 2520 490}%
\special{pa 2584 462}%
\special{pa 2561 456}%
\special{pa 2558 432}%
\special{pa 2520 490}%
\special{fp}%
\put(24.2000,-2.4000){\makebox(0,0)[lb]{$\pi$}}%
\end{picture}%
\label{SymmetryPIII1}
\caption{The transformations described in Theorem \ref{th:1} satisfy the relations:
${s_0}^2={s_1}^2={s_2}^2={\pi}^2=(s_0s_2)^2=(s_0s_1)^4=(s_1s_2)^4=1, \ {\pi}(s_0,s_1,s_2)=(s_2,s_1,s_0){\pi}.$}
\end{figure}
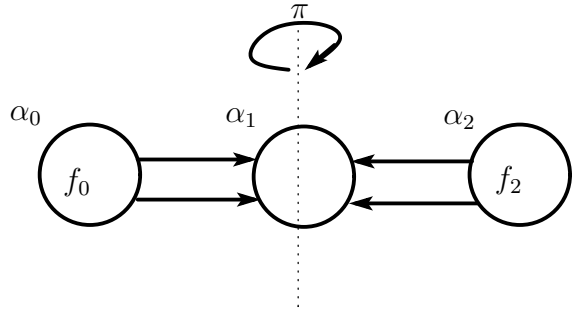

\begin{theorem}
For the system \eqref{symmetricPV} of type $A_3^{(1)}$, we make the change of parameters
\begin{equation}
\alpha_0=\beta_0, \ \alpha_1=-\frac{\eta}{a}, \ \alpha_2=\beta_2, \ \alpha_3=2\beta_1+\frac{\eta}{a}, \ \varphi=-2,
\end{equation}
from $\alpha_0,\alpha_1,\alpha_2,\alpha_3$ to $\beta_0,\beta_1,\beta_2$. This new system tends to the system \eqref{symmetricPIII} of type $C_2^{(1)}$ as $a \rightarrow 0$.
\end{theorem}
We note that
$$
\alpha_0+\alpha_1+\alpha_2+\alpha_3=\beta_0+2\beta_1+\beta_2=1.
$$

\end{document}